\documentclass[12pt]{article}

\usepackage{amsfonts,epsfig,graphicx}
\usepackage{amsmath,amssymb,amsthm}
\usepackage{epsf}

\newcommand{\F}{\mathcal{F}}
\newcommand{\R}{\mathcal{R}}

\newcommand{\SA}{\mathcal{S}}

\newcommand{\PA}{\mathcal{P}}
\newcommand{\X}{\mathcal{X}}
\newcommand{\G}{\mathcal{G}}
\newcommand{\M}{\mathcal{M}}
\newcommand{\1}{\textbf{1}}

\newcommand{\tE}{\tilde{E}}
\newcommand{\Q}{\mathcal{Q}}

\newcommand{\bSA}{{\bar{\SA}}}
\newcommand{\bU}{\bar{U}}
\newcommand{\bu}{\bar{u}}
\newcommand{\bx}{\bar{x}}
\newcommand{\by}{\bar{y}}
\newcommand{\bp}{\bar{p}}
\newcommand{\br}{\bar{r}}

\newcommand{\var}{\mbox{var}}

\usepackage{latexsym}

\usepackage[usenames,dvipsnames]{color}

\definecolor{britishracinggreen}{rgb}{0.0, 0.26, 0.15}

\newtheorem{thm}{Theorem}

\newtheorem{lem}{Lemma}
\newtheorem{cor}{Corollary}

\begin{document}
\begin{center}
\textbf{\Large A VARIATIONAL FORMULA FOR\\ \ \\ RISK-SENSITIVE REWARD}
\end{center}

\vspace{.5in}

\begin{center}
V.\ ANANTHARAM\footnote{EECS Department, University of California, Berkeley, CA 94720, USA. Research supported in part by the
ARO MURI grant W911NF- 08-1-0233, Tools for the
Analysis and Design of Complex Multi-Scale Networks, the
NSF grants CNS-0910702 and ECCS-1343398, and the NSF Science
\& Technology Center grant CCF-0939370, Science of Information.
A part of this work was done while this author was visiting IIT Bombay.} and V.\ S.\ BORKAR\footnote{Department of Elec.\ Engg., IIT Bombay, Powai, Mumbai 400076, India. Work supported in part by a J.\ C.\ Bose Fellowship and grant 11IRCCSG014 from IIT Bombay. A part of this work was done while this author was visiting the University of California, Berkeley.}
\end{center}

\vspace{1in}

\noindent \textbf{\large ABSTRACT:} We derive a variational formula for the optimal growth rate of reward in the infinite horizon risk-sensitive control problem for discrete time Markov decision processes with compact metric state and action spaces, extending a formula of
Donsker and Varadhan for the Perron-Frobenius eigenvalue of a positive operator.
This leads to a concave maximization formulation of the problem of determining this optimal growth rate.\\


\newpage

\section{Introduction}		\label{s:intro}

Infinite time horizon risk-sensitive control seeks to maximize the asymptotic growth rate for mean multiplicative reward in the standard Markov decision theory setting.
The optimal reward multiplier per step turns out to be the Perron-Frobenius eigenvalue of a positive $1$-homogeneous nonlinear operator.
The existence of this Perron-Frobenius eigenvalue
and an associated eigenfunction is ensured by the nonlinear Krein-Rutman theorem of
\cite[Theorem 3.1.1 and Proposition 3.1.5]{Ogiwara}
under suitable conditions (see also
\cite{Nussbaum}, \cite{MathSciNet}, \cite{Rajesh}, \cite{Chang}, \cite{Ari}).
Our aim here is to build on this nonlinear Krein-Rutman theorem to provide a variational formula for the optimal growth rate of reward in the spirit of the Donsker-Varadhan formula for the Perron-Frobenius eigenvalue of a nonnegative matrix \cite[section 3.1.2]{DemboZ},
\cite{DonVar}, \cite{Fried}.\\

Risk-sensitive control has traditionally been studied in the framework of cost minimization, 
see e.g. \cite{DiMasi}, \cite{Anna1}, \cite{Anna} for recent work on general state space models and \cite{Fleming1}, \cite{Hern} for its discrete state space precursors. Work on risk-sensitive reward maximization has been relatively uncommon, see e.g. \cite{Anna2}.
Unlike in the case of the classical discounted or ergodic costs, the two risk-sensitive control problems are not trivially equivalent by treating cost as a negative reward. In fact, risk-sensitive reward maximization is the natural set-up in portfolio optimization, see e.g. \cite{Cover}. Nevertheless, it has been commonplace to replace it by risk-sensitive cost minimization so as to exploit the vastly more abundant available machinery for the latter problem, see, e.g. equation (18) of \cite{Bielecki}. Interestingly, our approach is tailored for the risk-sensitive reward maximization problem. \\

The paper is organized as follows. This section presents the basic notation and
control-theoretic framework. In section \ref{s:principeigen} we develop
the role of the nonlinear
Krein-Rutman theorem in giving an expression for the optimal reward multiplier per
stage. In section \ref{s:variational} this is parlayed into a
variational expression for the optimal growth rate of reward. Theorem \ref{thm:mainstrong}
in section \ref{s:variational} is the main result of this paper.
Alternative variational formulations derived from the primary one are discussed in
section \ref{s:remarks}; each of these provides a different kind of insight into how
to think about the optimal growth rate of reward.
Some examples are worked out in section \ref{s:examples}
to illustrate the nature of the results. We close the paper with some concluding remarks in
section \ref{s:concremarks}.\\

We turn next to introducing our notation and the control-theoretic framework.
For a compact metric
space $\X$, $\M(\X)$ and $\PA(\X)$ will denote respectively the space of finite (signed) Borel measures on $\X$ and the space of probability measures on $\X$, both with the topology of weak convergence \cite{Billingsley}. $C(\X)$ will denote the Banach space of continuous maps $\X \mapsto \R$  with the supremum norm, denoted by
$\| \cdot \|$.
Thus $\M(\X)$ is the dual Banach space of $C(\X)$, with the weak-* topology
\cite{Rudin}. Let $\SA$ be a prescribed compact metric space called the \textit{state space} and $U$ another compact metric space, called the \textit{action space}. We shall consider an $\SA$-valued controlled Markov process $(X_n, n \geq 0)$ controlled by a $U$-valued control process $(Z_n, n \geq 0)$ defined as follows. Consider a complete probability space $(\Omega, \F, P)$ where $\Omega := (\SA\times U)^{\infty}$, and $\F$ is its product Borel $\sigma$-field. For $\omega = [(\omega_0, \omega'_0), (\omega_1, \omega'_1), (\omega_2, \omega'_2), \cdots] \in \Omega$ with $\omega_i \in \SA$ and $\omega'_i \in U \ \forall i$, define `canonical' random variables $X_i = \omega_i, Z_i = \omega'_i, i \geq 0$.  The probability measure $P$ on $(\Omega, \F)$ is then the law of $((X_n, Z_n), n \ge 0)$ defined as follows. The law of $X_0$ is prescribed and the law of $((X_n, Z_n), n \ge 0)$ is constructed inductively. For this purpose, define two increasing families of sub-$\sigma$-fields of $\F$ : $\F_n^- := \sigma(X_m, m \leq n; Z_m, m < n)$ and $\F_n := \sigma(X_m, m \leq n; Z_m, m \leq n)$  for $n \geq 0$.
First define the conditional law of $Z_0$ given $\F^-_0$ as
$\phi_0(du|X_0)$,
where
\begin{displaymath}
\phi_0(du|x_0) :\SA \mapsto \PA(U)
\end{displaymath}
is a prescribed kernel, i.e. $\phi_0(du|x)$ is a probability distribution in
$\PA(U)$ for all $x$ and $\phi_0(A|x)$ is Borel measurable in $x$
for all Borel subsets $A \subset U$.
Let $P_n$ denote the law of
$((X_0, Z_0), (X_1, Z_1), \cdots, (X_n, Z_n))$, defined as a probability measure on $(\Omega, \F_n)$, starting with $n = 0$. Define the law of $X_{n+1}$ given $\mathcal{F}_n$ as $p(dy|X_n, Z_n)$ where
\begin{displaymath}
p(dy | x, u) : \SA\times U \mapsto \PA(\SA)
\end{displaymath}
is a prescribed kernel, i.e. $p(dy|x,u)$ is a probability distribution
in $\PA(\SA)$ for all $(x,u) \in \SA \times U$ and
$p(A|x,u)$ is Borel measurable in $(x,u)$
for all Borel subsets $A \subset \SA$.
Define the conditional law of $Z_{n+1}$ given $\F^-_{n+1}$ as
\begin{displaymath}\phi_{n+1}(du|(X_0,Z_0), \cdots, (X_n,Z_n), X_{n+1})
\end{displaymath}
where
\begin{displaymath}
\phi_{n+1}(du|(x_0,u_0) \cdots, (x_n,u_n),x_{n+1}) : (\SA\times U)^n\times\SA \mapsto
\PA(U)
\end{displaymath}
is a prescribed kernel for each $n$. These together define $P_{n+1}$. By the Ionescu-Tulcea theorem (p.\ 101, \cite{Pollard}), we define a unique $P$ on $(\Omega, \F)$. By construction, for all Borel $A \subset \SA$,
\begin{eqnarray}
P(X_{n + 1} \in A | \F_n) &=& P(X_{n + 1} \in A | X_n, Z_n) \nonumber \\
&=& p(A | X_n, Z_n). \label{transition}
\end{eqnarray}
The $(Z_n, n \ge 0)$ constructed above will be referred to as admissible controls. We shall  also consider two special classes of admissible controls: \textit{stationary Markov controls} of the form
\begin{displaymath}
Z_n = v(X_n) \ \forall \ n,
\end{displaymath}
for some measurable $v : \SA \mapsto U$, and \textit{randomized stationary Markov controls} satisfying
\begin{displaymath}
P(Z_n \in A | \F_n) = P(Z_n \in A | X_n) = \varphi(A | X_n) \ \forall \ n, \forall \ \mbox{Borel} \ A \subset U,
\end{displaymath}
for some kernel $\varphi(du|x) : \SA \mapsto \PA(U)$. By a standard abuse of terminology, we identify these with the maps $v(\cdot), \varphi(\cdot | \cdot)$ resp. The sets thereof will be denoted by $SM$ and $RM$ respectively.  We view $SM$ as a subset of $RM$ by identifying $v(\cdot)$ with $\delta_{v(\cdot)}$, the Dirac measure at $v(\cdot)$.\\

The infinite horizon risk-sensitive reward we seek to characterize is
\begin{equation}		\label{growthrate}
\lambda :=  \sup_{x \in \SA} \sup \liminf_{N\uparrow\infty}\frac{1}{N}\log
E\left[e^{\sum_{m = 0}^{N - 1}r(X_m, Z_m, X_{m + 1})} | X_0 = x\right]~,
\end{equation}
where the second supremum is over all admissible controls.
Here $r(x,u,y)$ is
an extended-real-valued function on $\SA \times U \times \SA$,
called the  `per stage reward' on transitioning from $x$ to $y$ under action $u$. It should be noted that we will allow $e^{r(x,u,y)} = 0$ for some $(x,u,y)$, so
$r(x,u,y)$ should be thought
of as being allowed to take the extended real value $-\infty$.

Throughout the paper, we make the following assumptions about
$r(x,u,y)$ and $p(dy|x, u)$. We will occasionally explicitly recall these assumptions
to remind the reader of this.\\

\noindent \textbf{(A0)}: $e^{r(x,u,y)} \in C(\SA \times U \times \SA)$.\\

\noindent \textbf{(A1)}:
The maps $(x,u) \mapsto \int f(y)p(dy|x, u), f \in C(\SA)$ with $\|f\| \leq 1$, are
equicontinuous. This is true, e.g.,
if $\SA$ is a compact metric space,
$U$ is a compact metric space, and $p(dy|x,u) = \psi(y|x,u)\varphi(dy)$ with $\varphi \in \PA(\SA)$ having full support and $\psi(y| \cdot, \cdot), y \in \SA,$ equicontinuous.\\


We shall denote by $e^{r_M}$ the least upper bound for $e^{r(\cdot, \cdot, \cdot)}$, which
is finite by virtue of assumption \textbf{(A0)}.\\

Towards the end of the next section we will
build up to the main variational formula by first considering the
case where we have additional restrictions captured by the following assumptions.\\

\noindent \textbf{(A0+)}: Condition \textbf{(A0)} holds and we have
$e^{r(x,u,y)} > 0$ for all $(x,u,y)$.\\

\noindent \textbf{(A1+)}:
Condition \textbf{(A1)} holds and $p(dy | x, u)$ has full support for all $x, u$.
For instance,
if $\SA$ is a compact metric space,
$U$ is a compact metric space, and
$p(dy|x,u) = \psi(y|x,u)\varphi(dy)$ as above with $\psi(y| \cdot, \cdot), y \in \SA,$ equicontinuous, then
$\psi(\cdot | x,u) > 0$ on $\SA$ will ensure that this assumption holds. \\


We shall denote by  $e^{r_m} > 0$ the greatest lower bound for $e^{r(\cdot, \cdot, \cdot)}$
when \textbf{(A0+)} holds.

If $p(dx)$ and $q(dx)$ are finite nonnegative Borel measures on a compact metric space
$\X$, we write $D(p(dx) \| q(dx))$ for the relative entropy of $p(dx)$
with respect to $q(dx)$, defined by
\[
D(p(dx) \| q(dx)) := \begin{cases}
\int p(dx) \log l(x) & \mbox{ if we can write $p(dx) = l(x) q(dx)$}\\
\infty & \mbox{otherwise.}
\end{cases}
\]
See e.g. \cite{Van} for some of the basic properties of relative entropy.

\section{The Perron-Frobenius eigenvalue}		\label{s:principeigen}

Let assumptions \textbf{(A0)} and \textbf{(A1)} be in force.
 Define the operator $T : C(\SA) \mapsto C(\SA)$ by
\begin{equation}
Tf(x) := \sup_{\phi \in \PA(U)}\int\int p(dy | x, u)\phi(du)e^{r(x, u, y)}f(y). \label{operator}
\end{equation}
For fixed $x \in \SA$ on the left hand side of (\ref{operator})
the supremum on the right hand side is the expectation of a continuous affine function on a compact set of probability measures. Hence, it is a maximum attained at a Dirac measure.
For each fixed $f \in C(\SA)$,
a standard measurable selection theorem \cite[Lemma 1, p.\ 182]{Benes} allows us to choose
the family of maximizers, parametrized by $x \in \SA$, as a measurable function
$v : \SA \mapsto U$.
To see that $T$ is a map $C(\SA) \mapsto C(\SA)$, note that for $f \in C(\SA)$ with $\|f\| \leq R$,
\begin{eqnarray*}
\lefteqn{|Tf(x) - Tf(x')|} \\
 &=& |\sup_{\phi \in \PA(U)}\int\int p(dy | x, u)\phi(du)e^{r(x, u, y)}f(y) \\
 && - \ \sup_{\phi \in \PA(U)}\int\int p(dy | x', u)\phi(du)e^{r(x', u, y)}f(y)| \\
&=& |\sup_u\int p(dy | x, u)e^{r(x, u, y)}f(y) \\
&& - \ \sup_u\int p(dy | x', u)e^{r(x', u, y)}f(y)| \\
&\leq& e^{r_M}\sup_u\sup_{f: \|f\| \leq R}|\int p(dy | x, u)f(y) \\
&& - \ \int p(dy | x', u)f(y)| + R\max_{u,y}\left|e^{r(x, u, y)}- e^{r(x', u, y)}\right|.
\end{eqnarray*}
As $x \rightarrow x'$, the first term on the right tends to zero by \textbf{(A1)} and the second term on the right tends to zero by uniform continuity of
$e^r$,
being a continuous function defined on a compact set, by
\textbf{(A0)}.  In fact, this shows that $Tf, \|f\| \leq R$, are equicontinuous and bounded. Also, from  the definition of $T$, it is straightforward to check that
 \begin{displaymath}
 \|Tf - Tg\| \leq e^{r_M}\|f - g\|.
  \end{displaymath}
which establishes $T$ as a continuous (in fact, Lipschitz) map $C(\SA) \mapsto C(\SA)$.\\

  Likewise, define, for $f \in C(\SA)$,
  \begin{displaymath}
  T^{(n)}f(x) := \sup E\left[e^{\sum_{m = 0}^{n - 1}r(X_m, Z_m, X_{m + 1})}f(X_n) | X_0 = x\right],
   \end{displaymath}
where the supremum
is over all admissible control processes. Then $T^{(1)} = T$, by virtue of the
measurable selection theorem alluded to after (\ref{operator}).
We use the convention $T^{(0)} :=$ the identity map. \\

\begin{lem}			\label{lem:semigroup}
$(T^{(n)}, n \geq 0)$ is a semigroup of operators on $C(\SA)$.
\hfill $\Box$
\end{lem}

\noindent
\textit{Proof of Lemma \ref{lem:semigroup}}:
Note that we need to verify that $T^{(n)}$ for $n \ge 2$ maps $C(\SA)$ to
$C(\SA)$ as part of the stated claim. This follows as a corollary of the proof,
which establishes that $T^{(n)}$ is the $n$-fold concatenation of $T$ with itself.
The  proof follows by a standard dynamic programming argument. Specifically,
we first have
 \begin{eqnarray}
 &&T^{(n)}f(x) \nonumber \\
  &=& \sup E\left[e^{\sum_{m = 0}^{n - 1}r(X_m, Z_m, X_{m + 1})}f(X_n) | X_0 = x\right] \nonumber \\
 &\leq& \sup E\left[e^{r(X_0, Z_0, X_{1})}\sup E\left[e^{\sum_{m = 1}^{n - 1}r(X_m, Z_m, X_{m + 1})}f(X_n)|X_0, Z_0,X_1\right] | X_0 = x\right] \nonumber \\
 &=& \sup E\left[e^{r(X_0, Z_0, X_{1})}T^{(n-1)}f(X_1) | X_0 = x\right] , \label{one-1}
 \end{eqnarray}
 where the inner supremum in the second line is over the control sequence from time 1 onwards, conditioned on $X_0 = x_0, Z_0 = z_0, X_1 = x_1$ (say).
Secondly, let $\epsilon > 0$.
By \cite[Lemma 1, p.\ 55]{diffusion},
conditioned on $(X_0, Z_0, X_1)$,
there exists  an
admissible state-control sequence  $(X'_m, Z'_m), m \geq 1$, with
$X'_1 = X_1$ such that
 \begin{eqnarray*}
\lefteqn{E\left[e^{\sum_{m = 1}^{n - 1}r(X'_m, Z'_m, X'_{m + 1})}f(X'_n)|X'_1\right]}\\
&\geq&
\sup E\left[e^{\sum_{m = 1}^{n - 1}r(X_m, Z_m, X_{m + 1})}f(X_n)|X_1\right]
 - \epsilon, \ \mbox{a.s.}
 \end{eqnarray*}
 Let $X'_0 = X_0 = x, Z'_0 :=$ argmax$\left(\int p(dy|x,\cdot)e^{r(x,\cdot,y)}T^{(n-1)}f\right)$. Then
 \begin{displaymath}
 (X'_0, Z'_0), (X'_1, Z'_1), \cdots , (X'_n, Z'_n))
 \end{displaymath}
  is an admissible state-control sequence and

\begin{eqnarray}
 T^{(n)}f(x) &\geq& E\left[e^{\sum_{m = 0}^{n - 1}r(X'_m, Z'_m, X'_{m + 1})}f(X'_n) | X'_0 = x\right] \nonumber \\
 &\geq& E\left[e^{r(X_0, Z_0, X_{1})}\sup E\left[e^{\sum_{m = 1}^{n - 1}r(X_m, Z_m, X_{m + 1})}f(X_n)|X_1\right]\right] - e^{r_M}\epsilon \nonumber \\
 &=& E\left[e^{r(X_0, Z_0, X_{1})}T^{(n-1)}f(X_1) | X_0 = x\right] - e^{r_M}\epsilon \nonumber \\
&=& T^{(1)}\left(T^{(n-1)}f \right)(x) - e^{r_M}\epsilon \label{two-2}
\end{eqnarray}

 Combining (\ref{one-1}), (\ref{two-2}) and using the fact that $\epsilon > 0$ was arbitrary, it follows that $T^{(n)} = T^{(1)}\circ T^{(n-1)}$.
 A similar argument shows that $T^{(n)}f = T^{(n-1)}\circ Tf$. \hfill $\Box$

\ \\

The semigroup $(T^{(n)}, n \geq 0)$,  is precisely the discrete time counterpart of the Nisio semigroup \cite{Nisio}. \\

Let
$C^+(\SA) := \{f \in C(\SA) : f(x) \geq 0\}$
denote the set of nonnegative functions in $C(\SA)$. Then
$C^+(\SA)$ is a {\em cone}, i.e. it is
closed under addition and scalar multiplication by nonnegative real numbers, and we have
$C^+(\SA) \cap (- C^+(\SA)) = \{ \theta\}$ where $\theta$ denotes
the constant function that is identically zero.
Thus
$C^+(\SA)$ defines a partial order on $C(\SA)$, denoted $\ge$,
given by $f \ge g$ if
$f - g \in C^+(\SA)$.
We write $f > g$ (equivalently, $g < f$) if $f \ge g, f \neq g$, and
we write $f >> g$ if $f - g$ is a strictly positive function in $C(\SA)$ or equivalently if
$f - g \in \mbox{int}(C^+(\SA))$, where $\mbox{int}(C^+(\SA))$ denotes the
interior of $C^+(\SA)$.
The {\em dual cone} of
$C^+(\SA)$ is the cone in the
dual Banach space $\M(\SA)$ given by
$\{\mu \in \M(\SA) : \int fd\mu \geq 0 \ \forall \ f \in C^+(\SA)\}$.
This is the set of finite nonnegative measures on $\SA$, which we denote by
$\M^+(\SA)$.
For more on cones in Banach spaces, see \cite{Aliprantis}.


Let us now make the additional assumption \textbf{(A0+)} and  \textbf{(A1+)}.
One can then verify the following additional properties of $T^{(n)}$ for each $n \ge 1$.

\begin{enumerate}

\item $T^{(n)}$ is strictly increasing, i.e., $f < g$ implies $T^{(n)}f < T^{(n)}g$.
In view of the fact established above that $(T^{(n)}, n \geq 0)$
is a semigroup, it suffices to prove this claim for $n=1$.
We know that there is a measurable function
$v : \SA \mapsto U$ such that
\begin{displaymath}
Tf(x) = \int p(dy|x,v(x)) e^{r(x,v(x),y)}f(y)~.
\end{displaymath}
Then
\begin{eqnarray*}
\lefteqn{Tg(x) - Tf(x)} \\
&\geq& \int p(dy|x,v(x)) e^{r(x,v(x),y)}g(y) - \int p(dy|x,v(x)) e^{r(x,v(x),y)}f(y)\\
&\geq&e^{nr_m} \int p(dy|x,v(x)) (g(y) - f(y))\\
&>& 0~,
\end{eqnarray*}
because $f < g, f \neq g$ and support$(p(dy | x, u)) = \SA \ \forall \ x, u$.

\item $T^{(n)}$ is strongly positive, i.e., $f \in C^+(\SA), \ f \neq \theta \Longrightarrow T^{(n)}f \in$ int$(C^+(\SA))$. This follows from the fact that for any $u_0 \in U$,
    \begin{displaymath}
    T^{(n)}f(x) \geq e^{nr_m}\int p(dy|x, u_0)f(y) > 0,
    \end{displaymath}
    where we use the fact that support$(p(dy|x, u_0)) = \SA$.

\item $T^{(n)}$ is positively one-homogeneous, i.e., for $c > 0$, $T^{(n)}(cf) = cT^{(n)}f$.
(This holds under the weaker assumptions \textbf{(A0)} and  \textbf{(A1)}.)

\item For $M > e^{-nr_m}$ and $\breve{f} \in C(\SA)$ defined by $\breve{f}(\cdot) \equiv 1$, $MT^{(n)}\breve{f} > \breve{f}$.

\item $T^{(n)}$ is compact.
(This holds under the weaker assumptions \textbf{(A0)} and  \textbf{(A1)}.)
It suffices to verify this for $n = 1$, the general case being then a consequence of the semigroup property.
    By \textbf{(A1)}, the family
    $x \mapsto F_f(x, u) := \int f(y)e^{r(x,u,y)}p(dy | x,u), u \in U, \|f\| \leq R$, is equicontinuous and bounded in $C(\SA)$-norm by $e^{r_M}R$. Hence it is
 relatively compact in $C(\SA)$ by the Arzela-Ascoli theorem. Let $\delta \in [0,1] \mapsto w_{\delta}(\cdot)$ denote its common modulus of continuity relative to a compatible metric $\kappa$ on $\SA$. Then  $T : C(\SA) \mapsto C(\SA)$ satisfies $\|Tf\| \leq e^{r_M}R$ for $\|f\| \leq R, f \in C(\SA)$, and,
 \begin{eqnarray*}
&& \sup_{x, y \in \SA, \kappa(x,y) < \delta}\|Tf(x) - Tf(y)\| \\
&\leq& \sup_{x, y \in \SA, \kappa(x,y) < \delta}\|\sup_uF_f(x,u) - \sup_uF_f(y,u)\| \\
&\leq& \sup_{x, y \in \SA, \kappa(x,y) < \delta}\sup_u\|F_f(x,u) - F_f(y,u)\| \\
&\leq& w_{\delta}(F_f) \stackrel{\delta\downarrow 0}{\rightarrow} 0
\end{eqnarray*}
uniformly in $f: \|f\| \leq R$. Thus $Tf, \|f\| \leq R$, ie equicontinuous. By Arzela-Ascoli theorem, it is relatively compact, implying that $T: C(\SA) \mapsto C(\SA)$ is a compact operator.
\end{enumerate}


The preceding considerations allow us to state the following theorem.\\

\begin{thm}			\label{thm:principeigen}
Under the assumptions
\textbf{(A0+)} and  \textbf{(A1+)},
there exists a unique $\rho > 0$ (the \textit{Perron-Frobenius eigenvalue}) and a
$\psi \in$ int$(C^+(\SA))$ such  that $T\psi = \rho \psi$, i.e.,
\begin{equation}
\rho\psi(x) = \sup_{\phi \in \PA(U)}\int\int p(dy | x, u) \phi(du)e^{r(x, u, y)}\psi(y), \label{DP}
\end{equation}
with $\rho$ given by
\begin{eqnarray}
\rho &=&
\inf_{f \in {\rm int}(C^+(\SA)) }\sup_{\mu \in \M^+(\SA)}\frac{\int Tfd\mu}{\int fd\mu} \nonumber\\
&=&
\sup_{f \in {\rm int}(C^+(\SA)) }\inf_{\mu \in \M^+(\SA)}\frac{\int Tfd\mu}{\int fd\mu}~. \label{changformula}
\end{eqnarray}
\hfill $\Box$
\end{thm}

Equation (\ref{changformula}) is an abstract version of the celebrated Collatz-Wielandt formula for the Perron-Frobenius eigenvalue of irreducible nonnegative matrices,
see e.g. \cite{Meyer}.
\\

Before proceeding to the proof of Theorem \ref{thm:principeigen}, it is
appropriate to make a few remarks. A great deal is known about analogs of the Perron-Frobenius theorem for increasing positively one-homogeneous
maps on finite dimensional vector spaces, see the recent book \cite{LemmNuss}.
When the map is on an ordered Banach space (and one talks about a Krein-Rutman
theorem rather than a Perron-Frobenius theorem, in view of the seminal work
in \cite{KR}), we rely on
Theorem 3.1.1, Proposition 3.1.5, and Lemma 3.1.7 of \cite{Ogiwara}, as seen in the proof below
(see also \cite{Nussbaum}, \cite{MathSciNet}).
These results in \cite{Ogiwara} are themselves stated in a
much broader context than the special case of the Banach space $C(\SA)$
and the order structure defined by the cone
$C^+(\SA)$, with
$\SA$ a compact metric space,  which suffices for our purposes. The recent papers
\cite{Rajesh} and \cite{Chang} claim even stronger nonlinear Krein-Rutman theorems.
However, it has been recognized in \cite{Ari} that some of the claims in these
papers are wrong. The proof of the Theorem \ref{thm:principeigen} given
below does not rely in any way on \cite{Rajesh}, \cite{Chang}, or \cite{Ari}.\\

\textit{Proof of Theorem \ref{thm:principeigen}}: We define
\begin{displaymath}
\| T^{(n)} \|_+ := \sup \{ \| T^{(n)} f\| ~:~  f\in C^+(\SA), \| f\| \le 1 \}~,~~n \ge 0~.
\end{displaymath}
Since $(T^{(n)}, n \geq 0)$ is a positive semigroup, it is straightforward to check that
$\| T^{(k+l)} \|_+ \le \| T^{(k)} \|_+ \| T^{(l)} \|_+$ for all $k, l \ge 0$, and so
\begin{displaymath}
r(T) := \lim_{n \to \infty} \| T^{(n)} \|_+^{\frac{1}{n}}
\end{displaymath}
exists. By the fourth of the properties of the semigroup $(T^{(n)}, n \geq 0)$
shown above, we have $r(T) > 0$. It will turn out that the $\rho$ promised in
the statement of Theorem \ref{thm:principeigen} is just $r(T)$.

Strong positivity of $T$, which was shown above, verifies assumption A4 in
\cite[pg. 47]{Ogiwara}, and the facts that $T$ is compact (as established above), one-homogeneous, and order preserving are respectively the conditions A1, A2, and A3
in \cite[pg.47]{Ogiwara}. Thus \cite[Proposition 3.1.5.]{Ogiwara} provides the
additional requirement in the statement of \cite[Theorem 3.1.1]{Ogiwara} that $T$
have an eigenvalue, and \cite[Theorem 3.1.1]{Ogiwara} states that with
$\rho$ taken to be $r(T)$ there exists a $\psi \in \mbox{int}(C^+(\SA))$ such that
(\ref{DP}) holds.

It remains to establish (\ref{changformula}), where we now know that $\rho = r(T)$.
We have
\begin{displaymath}
\rho \ge
\inf_{f \in {\rm int}(C^+(\SA)) }\sup_{\mu \in \M^+(\SA)}\frac{\int Tfd\mu}{\int fd\mu}~,
\end{displaymath}
which comes from substituting $\psi$ as a choice for $f$ on the right hand
side. Similarly, we have
\begin{displaymath}
\rho \le
\sup_{f \in {\rm int}(C^+(\SA)) }\inf_{\mu \in \M^+(\SA)}\frac{\int Tfd\mu}{\int fd\mu}~.
\end{displaymath}
Thus it suffices to establish
\begin{equation}
\inf_{f \in {\rm int}(C^+(\SA)) }\sup_{\mu \in \M^+(\SA)}\frac{\int Tfd\mu}{\int fd\mu}
\ge \rho
\ge \sup_{f \in {\rm int}(C^+(\SA)) }\inf_{\mu \in \M^+(\SA)}\frac{\int Tfd\mu}{\int fd\mu}~.	\label{Desired}
\end{equation}
Given $f \in {\rm int}(C^+(\SA))$,
we have
\begin{displaymath}
Tf \le \left( \sup_{\mu \in \M^+(\SA)}\frac{\int Tfd\mu}{\int fd\mu} \right) f~.
\end{displaymath}
From \cite[Lemma 3.1.7 (ii)]{Ogiwara}, we have
$r(T) \le\sup_{\mu \in \M^+(\SA)}\frac{\int Tfd\mu}{\int fd\mu}$.
Since this holds for all $f \in {\rm int}(C^+(\SA))$,
this establishes the first inequality in (\ref{Desired}). The proof of the second inequality
in (\ref{Desired}) is similar, based on \cite[Lemma 3.1.7 (iii)]{Ogiwara}.
This concludes the proof of Theorem \ref{thm:principeigen}.
\hfill $\Box$

\ \\

Next we show that $\log\rho$ is in fact the optimal growth rate of the
risk-sensitive reward. For a development of the analogous result in the
case of controlled diffusion processes, see \cite{ABK}.
As argued earlier, in connection with the right hand side of (\ref{operator}),
for each $x \in \SA$,
the supremum on the right hand side of (\ref{DP})
is the expectation of a continuous affine function on a compact set of probability
measures, and is therefore a maximum attained at a Dirac measure.
A standard measurable selection theorem \cite[Lemma 1, p.\ 182]{Benes}
then allows us to identify
the family of maximizers, parametrized by $x \in \SA$,
with an element of $SM$, which we denote by $v^*(\cdot)$.
Letting $(X_n^*, n \ge 0)$ denote the chain governed by the stationary
Markov strategy $v^*(\cdot)$ and $(Z^*_n = v^*(X^*_n), n \ge 0)$ the corresponding control sequence,
we then have
\begin{displaymath}
\rho \psi(x) = E\left[e^{r(x, v^*(x), X_1^*)}\psi(X_1^*)\right],
\end{displaymath}
and, more generally, by iterating, we have, for all $x \in \SA$,
\begin{displaymath}
\rho^n\psi(x) =  E\left[e^{\sum_{m=0}^{n - 1}r(X_m^*, Z_m^*, X_{m + 1}^*)}\psi(X_n^*)| X_0^* = x\right].
\end{displaymath}

Since $\psi(x) \in \mbox{int}(C^+(\SA))$,
we have $0 < c < \psi(\cdot) < C < \infty$ for some constants $c, C$ when $\psi$ is chosen with, say, $\|\psi\| = 1$. Thus, for all $x \in \SA$,
\begin{displaymath}
\frac{c}{C}E\left[e^{\sum_{m=0}^{n - 1}r(X_m^*, Z_m^*, X^*_{m + 1})} | X_0^* = x\right] \leq \rho^n \leq \frac{C}{c}E\left[e^{\sum_{m=0}^{n - 1}r(X_m^*, Z_m^*, X^*_{m + 1})} | X_0^* = x\right].
\end{displaymath}
Hence
\begin{displaymath}
\log\rho = \lim_{n\uparrow\infty}\frac{1}{n}\log E\left[e^{\sum_{m = 0}^{n - 1}r(X^*_m, Z_m^*, X^*_{m + 1})} | X_0^* = x\right].
\end{displaymath}
For any other admissible state-control sequence $((X_n, Z_n), n \ge 0)$, we have
\begin{displaymath}
\rho\psi(x) \leq E\left[e^{r(x, Z_0, X_1)}\psi(X_1) | X_0 = x\right].
\end{displaymath}
Iterating,
\begin{displaymath}
\rho^n\psi(x) \leq E\left[e^{\sum_{m=0}^{n - 1}r(X_m, Z_m, X_{m + 1})}\psi(X_n) | X_0 = x\right].
\end{displaymath}
 and therefore
\begin{displaymath}
\log\rho \leq \liminf_{n\uparrow\infty}\frac{1}{n}\log E\left[e^{\sum_{m = 0}^{n - 1}r(X_m, Z_m, X_{m + 1})} | X_0 = x\right].
\end{displaymath}
We have proved:\\

\begin{thm}		\label{thm:growthrate}
Under the assumptions
\textbf{(A0+)} and  \textbf{(A1+)}, we have, for all $x \in \SA$,
\begin{displaymath}
\log\rho = \sup\liminf_{n\uparrow\infty}\frac{1}{n}\log E\left[e^{\sum_{m = 0}^{n - 1}r(X_m, Z_m, X_{m + 1})} | X_0 = x\right],
\end{displaymath}
where the supremum on the right is over all admissible controls and
$\rho$ on the left is given as in Theorem \ref{thm:principeigen}. Furthermore, this supremum is a maximum attained at some $v^*(\cdot) \in SM$.
\hfill $\Box$
\end{thm}

An immediate consequence is the following.

\begin{cor}			\label{cor:growthrate}
Under the assumptions
\textbf{(A0+)} and  \textbf{(A1+)} we have
\begin{displaymath}
\lambda = \log \rho~,
\end{displaymath}
where $\lambda$ is the optimal growth rate of reward, as defined in
(\ref{growthrate}), and $\rho$ is as
defined in Theorem \ref{thm:principeigen}.
\hfill $\Box$
\end{cor}

\section{A variational formula}		\label{s:variational}

By (\ref{changformula}), we have
\begin{eqnarray*}
\rho &=& \inf_{f >> 0}\sup_{\mu \in \mathcal{M}^+(S):\int fd\mu = 1}\int\mu(dx)\sup_u\int p(dy|x,u)e^{r(x,u,y)}f(y) \\
 &=& \inf_{f >> 0}\sup_{\nu \in \mathcal{P}(S)}\int\nu(dx)\left(\frac{\sup_u\int p(dy|x,u)e^{r(x,u,y)}f(y)}{f(x)}\right) \\
 &=& \inf_{f >> 0}\sup_x\left(\frac{\sup_u\int p(dy|x,u)e^{r(x,u,y)}f(y)}{f(x)}\right) \\
 &=& \inf_{f >> 0}\sup_x\sup_u\int p(dy|x,u)e^{r(x,u,y) + \log f(y) - \log f(x)} \\
 &=& \inf_{f >> 0}\sup_{\gamma \in \PA(\SA \times U)}\int\int\gamma(dx,du)\int p(dy|x,u)e^{r(x,u,y) + \log f(y) - \log f(x)}.
 \end{eqnarray*}
 Introduce the notation
 \begin{eqnarray*}
\eta(dx,du,dy) &=& \eta_0(dx)\eta_1(du|x)\eta_2(dy|x,u) \\
&=& \tilde{\eta}(dx,du)\eta_2(dy|x,u).
\end{eqnarray*}
Let
\begin{eqnarray*}
\mathcal{G} &:=& \{\eta(dx,du,dy) : \eta_0 \ \mbox{is invariant under the transition kernel}\\
&& \int_U \eta_2(dy|x,u)\eta_1(du|x)\},
\end{eqnarray*}
i.e. $\eta \in \mathcal{G}$ iff
\begin{displaymath}
\int \tilde{\eta}(dx, du) \eta_2(dy|x,u) = \eta_0(dy)~.
\end{displaymath}
Recall that $D( \cdot || \cdot )$ is convex and lower semi-continuous in both arguments \cite{Van}. Then

 \begin{eqnarray*}
 \lefteqn{\log \rho} \\
&=& \inf_{f >> 0}\sup_{\gamma}\log\int\int\int\gamma(dx,du)p(dy|x,u)e^{r(x,u,y) + \log f(y) - \log f(x)} \\
&=& \inf_{g \in C(\SA)}\sup_{\gamma}\log\int\int\int\gamma(dx,du)p(dy|x,u)e^{r(x,u,y) + g(y) - g(x)} \\
&=& \inf_{g \in C(\SA)} \sup_{\gamma} \sup_{\eta}\int\int\int\eta(dx,du,dy)\Big(r(x,u,y) + g(y) - g(x) \Big)\\
&& - \ D(\eta(dx,du,dy)||\gamma(dx,du)p(dy|x,u)) \\
&& \ \  (\mbox{by the Gibbs variational formula (Prop.\ 1.4.2(a), pp.\ 33-34, \cite{Dup}}) \\
&=& \sup_{\gamma}\sup_{\eta}\inf_{g \in C(\SA)}\int\int\int\eta(dx,du,dy)\Big(r(x,u,y) + g(y) - g(x) \Big)\\
&& - \ D(\eta(dx,du,dy)||\gamma(dx,du)p(dy|x,u)) \\
&&  \ \ \cdots\cdots (\mbox{by the min-max theorem \cite{Fan}}) \\
&=& \sup_{\gamma}\sup_{\eta}\inf_{g \in C(\SA)}\Big(\int\int\int\eta(dx,du,dy)\Big(r(x,u,y) + g(y) - g(x)\Big) \\
&& - \Big(D(\tilde{\eta}(dx,du)||\gamma(dx,du)) +  \int\int\tilde{\eta}(dx,du)D(\eta_2(dy|x,u)||p(dy|x,u)))\Big) \\
&=& \sup_{\eta}\inf_{g \in C(\SA)}\Big(\int\int\int\eta(dx,du,dy)\Big(r(x,u,y) + g(y) - g(x)\Big) \\
&& - \  \int\int\tilde{\eta}(dx,du)D(\eta_2(dy|x,u)||p(dy|x,u)) \\
&& \ \ \cdots\cdots (\mbox{by setting} \ \gamma = \tilde{\eta}) \\
&=& \sup_{\eta \in \mathcal{G}}\Big[\inf_{g \in C(\SA)}\Big(\int\int\int\eta(dx,du,dy)\Big(r(x,u,y) + g(y) - g(x)\Big) \\
&& - \  \int\int\tilde{\eta}(dx,du)D(\eta_2(dy|x,u)||p(dy|x,u))\Big] \\
&& \ \ \cdots\cdots (\mbox{because} \ \Big[\cdots\Big] = -\infty \ \forall \ \eta \notin \mathcal{G}) \\
&=& \sup_{\eta \in \mathcal{G}}\Big(\int\int\int\eta(dx,du,dy)r(x,u,y)  \\
&& - \  \int\int\tilde{\eta}(dx,du)D(\eta_2(dy|x,u)||p(dy|x,u))\Big) \\
&& \ \ \cdots\cdots (\mbox{because} \ \eta \in \mathcal{G} \Longrightarrow \int\eta(dx,du,dy)(g(y) - g(x)) = 0) \\
\end{eqnarray*}
Thus we have:\\

\begin{thm}		\label{thm:mainweak}
Under the assumptions
\textbf{(A0+)} and  \textbf{(A1+)}, the optimal growth rate of reward $\lambda$,
as defined in (\ref{growthrate}), has the variational characterization
\begin{eqnarray}
\lambda = \log\rho &=& \sup_{\eta \in \mathcal{G}}\Big(\int\int\int\eta(dx,du,dy)r(x,u,y) \nonumber \\
&& - \  \int\int\tilde{\eta}(dx,du)D(\eta_2(dy|x,u)||p(dy|x,u))\Big)~, \label{varweak}
\end{eqnarray}
where $\rho$ is defined as in Theorem \ref{thm:principeigen}.
\hfill $\Box$
\end{thm}

The following result, which uses a limiting argument to strengthen Theorem
\ref{thm:mainweak}, is the main result of this paper.

\begin{thm}		\label{thm:mainstrong}
Under the assumptions
\textbf{(A0)} and  \textbf{(A1)}, the optimal growth rate of reward $\lambda$,
as defined in (\ref{growthrate}), has the variational characterization
\begin{eqnarray}
\lambda &=& \sup_{\eta \in \mathcal{G}}\Big(\int\int\int\eta(dx,du,dy)r(x,u,y) \nonumber \\
&& - \  \int\int\tilde{\eta}(dx,du)D(\eta_2(dy|x,u)||p(dy|x,u))\Big)~. \label{varstrong}
\end{eqnarray}
\hfill $\Box$
\end{thm}

Before proving Theorem \ref{thm:mainstrong}, let us first consider the
uncontrolled case. We can fit this into our framework by taking $U$ to be
a set with one point, so that
$p(dy|x,u) = \tilde{p}(dy|x)$ for all $u \in U$, for some kernel $\tilde{p}(dy|x)$,
and $r(x,u,y) = \tilde{r}(x,y)$ for all $u \in U$, for some $\tilde{r}(\cdot, \cdot)$.
Theorem \ref{thm:mainstrong} then specializes to the statement that the
growth rate of the reward, under the respective specializations of
conditions \textbf{(A0)} and  \textbf{(A1)},
is given by
\begin{eqnarray*}
\lambda
&=& \sup_{\alpha \in \tilde{\mathcal{G}}}\Big(\int\int\int\alpha(dx,dy) \tilde{r}(x,y)  \\
&& - \  \int\int\alpha_0(dx)D(\alpha_1(dy|x)||\tilde{p}(dy|x))\Big)
\end{eqnarray*}
where $\alpha(dx,dy) = \alpha_0(dx)\alpha_1(dy|x)$ and
\begin{displaymath}
\tilde{\mathcal{G}} := \{\alpha(dx,dy) = \alpha_0(dx)\alpha_1(dy|x) : \int\alpha_0(dx)\alpha_1(dy|x) = \alpha_0(dy)\}.
\end{displaymath}
This is then a version of the Donsker-Varadhan formula for the
Perron-Frobenius eigenvalue
of a positive operator \cite{DemboZ}, \cite{DonVar}, \cite{Fried}.
\\

\textit{Proof of Theorem \ref{thm:mainstrong}}: Let $\gamma(dy)$ be an arbitrary
probability distribution on $\SA$ with full support, and, for all $\epsilon > 0$ sufficiently
small, define the kernel
\[
p_\epsilon(dy|x,u) := \frac{1}{a(x,u) + \epsilon} \Big( e^{r(x,u,y)} p(dy|x,u) + \epsilon \gamma(dy) \Big)~,
\]
and the reward
\[
r_\epsilon(x,u,y) := \log (a(x,u) + \epsilon)~,
\]
where
\[
a(x,u) := \int e^{r(x,u,y)} p(dy|x,u)~.
\]
Since this kernel and reward satisfy the conditions
\textbf{(A0+)} and  \textbf{(A1+)}, we have from Theorem \ref{thm:mainweak}
that the optimal growth rate of reward for the risk-sensitive reward maximization
problem for this kernel and reward, call it $\lambda_\epsilon$, is given by
\begin{eqnarray}
\lambda_\epsilon &=&
\sup_{\eta \in \mathcal{G}}\Big(\int\int\int\eta(dx,du,dy)r_\epsilon(x,u,y) \nonumber \\
&& - \  \int\int\tilde{\eta}(dx,du)D(\eta_2(dy|x,u)||p_\epsilon(dy|x,u))\Big)~. \label{vareps}
\end{eqnarray}

From the formulation of the risk-sensitive objective we see
that $\lambda_\epsilon$ is nondecreasing in $\epsilon$, and that
$\lambda_\epsilon \ge \lambda$ for all $\epsilon > 0$,
where $\lambda$ is defined
as in (\ref{growthrate}).
This can be seen by writing the expression for the $n$-step
multiplicative
reward, i.e.
\[
E_\epsilon \left[e^{\sum_{m = 0}^{N - 1}r_\epsilon(X_m, Z_m, X_{m + 1})} | X_0 = x\right]~,
\]
as a multiple integral,
which reveals that this quantity is monotonically nondecreasing in $\epsilon$
for any initial condition $x \in \SA$ and any admissible control strategy.
Thus $\lim_{\epsilon \to 0} \lambda_\epsilon$ exists
and satisfies
\begin{equation}		\label{lambdalimit}
\lim_{\epsilon \to 0} \lambda_\epsilon \ge \lambda~.
\end{equation}

To prove (\ref{varstrong}), we will first prove that
\begin{eqnarray}
\lim_{\epsilon \to 0} \lambda_\epsilon &\le& \sup_{\eta \in \mathcal{G}}\Big(\int\int\int\eta(dx,du,dy)r(x,u,y) \nonumber \\
&& - \  \int\int\tilde{\eta}(dx,du)D(\eta_2(dy|x,u)||p(dy|x,u))\Big)~, \label{varleq}
\end{eqnarray}
and then prove that
\begin{eqnarray}
\lambda &\ge& \sup_{\eta \in \mathcal{G}}\Big(\int\int\int\eta(dx,du,dy)r(x,u,y) \nonumber \\
&& - \  \int\int\tilde{\eta}(dx,du)D(\eta_2(dy|x,u)||p(dy|x,u))\Big)~. \label{lambdalower}
\end{eqnarray}
Together with (\ref{lambdalimit}), these two claims establish (\ref{varstrong}).

For fixed $\eta \in \mathcal{G}$, let $\Psi_\epsilon(\eta)$
denote the expression inside the outer brackets on the right hand side
of (\ref{vareps}).
Then one has
\begin{equation}		\label{psieps}
\Psi_\epsilon(\eta) = - \  \int\int\tilde{\eta}(dx,du)
D(\eta_2(dy|x,u)||e^{r(x,u,y)} p(dy|x,u) + \epsilon \gamma(dy))~.
\end{equation}
Similarly, for fixed $\eta \in \mathcal{G}$, let $\Psi_0(\eta)$
denote the expression inside the outer brackets on the right hand side
of (\ref{varstrong}). We have
\begin{equation}		\label{psi}
\Psi_0(\eta) = - \  \int\int\tilde{\eta}(dx,du)
D(\eta_2(dy|x,u)||e^{r(x,u,y)} p(dy|x,u))~.
\end{equation}

In fact, (\ref{psieps}) reveals that
for each $\eta \in \mathcal{G}$ we have
$\Psi_\epsilon(\eta)$ is nondecreasing in $\epsilon$, and together with
(\ref{psi}), reveals that
for all $\epsilon > 0$ and $\eta \in \mathcal{G}$, we have
$\Psi_\epsilon(\eta) \ge \Psi_0(\eta)$. Thus we may conclude that for each
$\eta \in \mathcal{G}$ the limit $\lim_{\epsilon \to 0} \Psi_\epsilon(\eta)$ exists, and
that this limit satisfies $\lim_{\epsilon \to 0} \Psi_\epsilon(\eta) \ge \Psi_0(\eta)$.

Now, for all $\epsilon > 0$ and $\delta > 0$ sufficiently small, choose
$\eta^\delta_\epsilon \in \mathcal{G}$
such that $\Psi_\epsilon(\eta^\delta_\epsilon) > \lambda_\epsilon - \delta$.
Since $\mathcal{G}$ is compact, there is a decreasing sequence $(\epsilon_m, m \ge 1)$
with $\lim_{m \to \infty} \epsilon_m = 0$, such that the sequence
$(\eta^\delta_{\epsilon_m}, m \ge 1)$ has a limit in $\PA(\SA \times U \times \SA)$,
call it $\eta^\delta$. Further, since $\mathcal{G}$ is closed, we have
$\eta^\delta \in \mathcal{G}$. By the lower semicontinuity of
$D( \cdot \| \cdot)$ as a function of $(\cdot, \cdot)$ \cite{Van}
we have
\[
\sup_{\eta \in \mathcal{G}} \Psi_0(\eta) \ge \Psi_0(\eta^\delta)
\ge \lim_{m \to \infty} \Psi_{\epsilon_m}(\eta^\delta_{\epsilon_m})
\ge \lim_{m \to \infty} \lambda_{\epsilon_m} - \delta
= \lim_{\epsilon \to 0} \lambda_\epsilon - \delta~.
\]
Since $\delta > 0$ was arbitrary, this establishes (\ref{varleq}).

It remains to prove (\ref{lambdalower}). If $\sup_{\eta \in \mathcal{G}} \Psi_0(\eta)$
(i.e. the right hand side of (\ref{lambdalower})) equals $- \infty$ then there is nothing
to prove, so we may assume that this is not the case.
Given $\eta \in \mathcal{G}$ for which $\Psi_0(\eta) \neq - \infty$, consider
implementing the stationary Markov strategy defined by the kernel $\eta_1(du|x)$.
The expected multiplicative reward after $n$ steps when implementing this strategy,
conditioned on starting with the initial distribution $\eta_0(dx_0)$, is
\[
\int \cdot\cdot \int \eta_0(dx_0) \prod_{m=0}^{n-1}  \eta_1(du_m|x_m) p(dx_{m+1}| x_m, u_m)
e^{r(x_m,u_m,x_{m+1})}~.
\]
Since $\eta_2(dy|x,u)$ is absolutely continuous with respect to
$p(dy|x,u)$ for almost all $(x,u)$, this equals
\[
\int \cdot\cdot \int \eta_0(dx_0) \prod_{m=0}^{n-1}  \eta_1(du_m|x_m) \eta_2(dx_{m+1}| x_m, u_m)
e^{r(x_m,u_m,x_{m+1})}e^{ - \log \frac{\eta_2(d x_{m+1}|x_m,u_m)}
{p(d x_{m+1}|x_m,u_m)}}~.
\]
Let $\{X'_n\}$ denote a controlled Markov chain with controlled transition kernel $\eta_2(dy|x,u)$, initial law $\eta_0$, and controlled by $\eta_1(du|x) \in RM$. Then
\begin{eqnarray*}
\lefteqn{\lambda} \\
&\ge&~\lim_{n \to \infty} \frac{1}{n} \log \Big(\int \cdot\cdot \int \eta_0(dx_0) \prod_{m=0}^{n-1}  \eta_1(du_m|x_m) p(dx_{m+1}| x_m, u_m)
e^{r(x_m,u_m,x_{m+1})}\Big) \\
&\ge&~\lim_{n \to \infty} \frac{1}{n} \log \Big(
\int \cdot\cdot \int \eta_0(dx_0) \prod_{m=0}^{n-1}  \eta_1(du_m|x_m) \eta_2(dx_{m+1}| x_m, u_m)
\\
&& \times \ e^{r(x_m,u_m,x_{m+1}) - \log \frac{\eta_2(d x_{m+1}|x_m,u_m)}
{p(d x_{m+1}|x_m,u_m)}}
\Big)\\
&=&~\lim_{n \to \infty} \frac{1}{n} \log \Big(E\left[e^{\sum_{m=0}^{n-1}(r(X_m',Z_m',X_{m+1}') - \log\frac{d\eta_2(\cdot|X'_m,Z'_m)}
{dp(\cdot|X'_m,Z'_m)}(X_m'))}\right]\Big) \\
&\geq&~\lim_{n \to \infty} \frac{1}{n} E\left[\sum_{m=0}^{n-1}(r(X_m',Z_m',X_{m+1}') - \log\frac{d\eta_2(\cdot|X'_m,Z'_m)}
{dp(\cdot|X'_m,Z'_m)}(X_m'))\right] \\
&& \ \ \ \ \ \ \ \ \ \ \ \ \ \ \ \ \ \  (\mbox{by Jensen's inequality}) \\
&=& \Psi_0(\eta) \\
&& \ \ \ \ \ \ \ \ \ \ \ \ \ \ \ \ \ \  (\mbox{because} \ \eta \in \G).
\end{eqnarray*}

It follows that $\lambda$, as defined in (\ref{growthrate}), satisfies (\ref{lambdalower}),
which concludes the proof of Theorem \ref{thm:mainstrong}.
\hfill $\Box$.

\section{Remarks}			\label{s:remarks}

\begin{enumerate}

\item

Assume \textbf{(A0), (A1)}. Fix $\varphi \in RM$, and consider
$\{(X_n, Z_n), n \geq 0\}$ governed by the randomized stationary Markov strategy $\varphi$ as an uncontrolled $\SA\times U$-valued Markov chain. To be precise, let $\bSA$ denote
$\SA \times U$, let $\bU := \{ \bu \}$ be a one point set, and define
$\bp(d \by| \bx,\bu) : \bSA \times \bU \mapsto \PA(\bSA)$
by
\begin{displaymath}
\bp(d \by| \bx,\bu) := p(dy|x,u) \varphi(du^\prime|y)~,
\end{displaymath}
where $\bx := (x,u)$ and $\by := (y,u^\prime)$. Also, let
\begin{displaymath}
\br(\bx, \bu, \by) := r(x,u,y)~.
\end{displaymath}
It is straightforward to check that the assumptions \textbf{(A0), (A1)}
hold for the $\bSA$-valued chain with trivial control space $\bU$ and with the
transition kernel and one step reward as above.

Given
$\tau(dx, du, dy, du^\prime) = \tau_0(dx) \tau_1(du|x) \tau_2(dy|x,u)
\tau_3(du^\prime |x,u ,y)$, write $\tilde{\tau}(dx,du)$ for
$\tau_0(dx) \tau_1(du|x)$ and $\hat{\tau}(dy, du^\prime|x, u)$ for
$\tau_2(dy|x,u)
\tau_3(du^\prime |x,u ,y)$. Let
\begin{displaymath}
\mathcal{G}_+ := \{ \tau(dx,du,dy,du^\prime) ~:~
\int \int \tilde{\tau}(dx,du) \hat{\tau}(dy, du^\prime|x, u) =
\tilde{\tau}(dy, du^\prime) \}~.
\end{displaymath}

Further, given $\tau(dx,du,dy,du^\prime)$, we define
$\tau^\prime(dx,du,dy,du^\prime)$ by setting
\begin{displaymath}
\tau^\prime_0 :=\tau_0,~\tau^\prime_1 := \tau_1,~\tau^\prime_2 := \tau_2,~
\tau^\prime_3(d u^\prime|x,u,y) := \tau_1(du^\prime|y)~,
\end{displaymath}
with the corresponding definitions for $\tilde{\tau}', \hat{\tau}'$. We claim that $\tau^\prime \in \mathcal{G}_+$.
To see this, first observe that
$\int \int \tilde{\tau}(dx,du) \hat{\tau}(dy, du^\prime|x, u) =
\tilde{\tau}(dy, du^\prime)$ when integrated over $u^\prime$ gives
$\int \int \tilde{\tau}(dx,du) \tau_2(dy|x, u) =
\tau_0(dy)$. This means
\begin{eqnarray*}
\int \int \tilde{\tau}^\prime(dx,du) \hat{\tau}^\prime(dy, du^\prime|x, u) &=&
\int \int \tilde{\tau}(dx,du) \tau_2(dy|x, u) \tau_1(du^\prime|y)\\
&=& \tau_0(dy) \tau_1(du^\prime|y)\\
&=& \tilde{\tau}(dy,du^\prime) = \tilde{\tau}^\prime(dy,du^\prime)~,
\end{eqnarray*}
which establishes the claim.

 Let $\lambda_{\varphi}$ denote the asymptotic growth rate of reward under
the fixed randomized stationary Markov strategy $\varphi$.
Then by applying Theorem \ref{thm:mainstrong} to the
$\bSA$-valued chain with trivial control space $\bU$ defined above, we have
\begin{eqnarray}
\lambda_{\varphi} &=& \sup_{\tau \in \mathcal{G}_+}\Big(\int\int\int\tau(dx,du,dy,U)r(x,u,y) - \nonumber \\
&& \int\int\tilde{\tau}(dx,du)D(\hat{\tau}(dy,du'|x,u)||p(dy|x,u)\varphi(du'|y))\Big). \label{var2}
\end{eqnarray}
Then we have
\begin{eqnarray*}
\lefteqn{\sup_{\varphi}\lambda_{\varphi}} \\
&=&\sup_{\varphi \in RM}\sup_{\tau \in \mathcal{G}_+}\Big(\int\int\int\tau(dx,du,dy,U)r(x,u,y) - \\
&& \int\int\tilde{\tau}(dx,du)D(\tau_2(dy|x,u)\tau_3(du'|x,u,y)||p(dy|x,u)\varphi(du'|y))\Big) \\
&\stackrel{(a)}{=}& \sup_{\tau \in \mathcal{G}_+}\Big(\int\int\int\tau^\prime(dx,du,dy,U)r(x,u,y) - \\
&& \int\int\tilde{\tau}^\prime(dx,du)D(\tau^\prime_2(dy|x,u)||p(dy|x,u))\Big) \\
&\stackrel{(b)}{=}& \sup_{\eta \in \mathcal{G}}\Big(\int\int\int \eta(dx,du,dy)r(x,u,y) - \\
&& \int\int\tilde{\eta}(dx,du)D(\eta_2(dy|x,u)||p(dy|x,u))\Big)\\
&=& \lambda~.
\end{eqnarray*}
Here, to justify step (a), notice that for every $\tau \in \mathcal{G}_+$, we have
shown that $\tau^\prime \in \mathcal{G}_+$. Therefore we have both
\begin{displaymath}
\int\int\int\tau^\prime(dx,du,dy,U)r(x,u,y) =
\int\int\int\tau (dx,du,dy,U)r(x,u,y)
\end{displaymath}
and
\begin{eqnarray*}
&&\int\int\tilde{\tau}^\prime(dx,du)D(\tau^\prime_2(dy|x,u)||p(dy|x,u)) \\
&=& \int\int\tilde{\tau}(dx,du)D(\tau_2(dy|x,u)||p(dy|x,u)).
\end{eqnarray*}
The choice of
$\varphi(du^\prime|y) = \tau_1(du^\prime|y)$ (which also equals
$\tau^\prime_3(du^\prime|x,u,y)$) would make the expression
\begin{displaymath}
\int \int \int \tilde{\tau}^\prime(dx,du) \tau_2^\prime(dy|x,u)
D(\tau^\prime_3(du^\prime|x,u,y)||\varphi(du^\prime|y))
\end{displaymath}
equal to zero, whereas
the expression
\begin{displaymath}
\int \int \int \tilde{\tau}(dx,du) \tau_2(dy|x,u)
D(\tau_3(du^\prime|x,u,y)||\varphi(du^\prime|y))
\end{displaymath}
is nonnegative.
To justify step (b) note that for every $\tau \in \mathcal{G}_+$,
we have $\tau_0(dx)\tau_1(du|x)\tau_2(dy|x,u) \in \mathcal{G}$, and
conversely for every $\eta \in \mathcal{G}$ we get $\tau \in \mathcal{G}_+$
by defining $\tau(dx,du,dy,du^\prime) :=\eta(dx,du,dy)\eta_1(du^\prime|y)$.
Furthermore, this $\tau$ satisfies $\tau^\prime = \tau$.

The upshot is that we have proved
\begin{equation}
\lambda  = \sup_{\varphi \in RM} \lambda_\phi~. \label{eq:lambda}
\end{equation}

Under \textbf{(A0+), (A1+)}, this supremum is in fact a maximum by virtue of Theorem 2.\\

\item Since $D( \cdot || \cdot)$ is convex and lower semi-continuous in its arguments as noted earlier, (\ref{varstrong}) is a concave maximization problem on the convex\footnote{See \cite[section 11.2.3, p.\ 358]{Bork} for the proof of convexity} set\\
\begin{eqnarray*}
\G_1&:=& \{\eta(dx)\varphi(du|x)\mu(dy|x,u) : \eta \ \mbox{is invariant under the transition} \\
&& \mbox{kernel} \ x \mapsto \int_U\varphi(du|x)\mu(dy|x,u)\}.
\end{eqnarray*}
It is worthwhile to compare this formulation with the classical dynamic programming approach.
Recall that the dynamic programming equation (\ref{DP}) is the nonlinear eigenvalue problem
\begin{displaymath}
\rho V(x) = \sup_{\varphi}\left(\int\int p(dy|x, u)\varphi(dy|u)e^{r(x, u, y)}V(y)\right).
\end{displaymath}
Consider the standard `log transformation' $\zeta(x) := \log V$. Then
\begin{displaymath}
\log\rho + \zeta(x)
= \sup_{\varphi}\log\left(\int\int p(dy|x, u)\varphi(du|x)e^{r(x, u, y) + \zeta(y)}\right).
\end{displaymath}
We treat $x$ as a fixed parameter on the right hand side. By the Gibbs variational principle, we have
\begin{eqnarray}
\lefteqn{\log\rho + \zeta(x)} \nonumber \\
&=& \sup_{\varphi}\sup_{\mu(\cdot,\cdot|x) \in \PA(U\times\SA)}\Big(\int\mu(du,dy|x)(r(x, u, y) + \zeta(y)) -  \nonumber \\
&& D(\mu(du,dy|x)||p(dy|x,u)\varphi(du|x))\Big).  \label{DPnew}
\end{eqnarray}
Equation (\ref{DPnew}) is the dynamic programming equation for an ergodic team problem whose  `per stage payoff' function is
\begin{displaymath}
r(x, u, y) - D(\mu(du,dy|x)||p(dy|x,u)\varphi(du|x)),
\end{displaymath}
where $\mu$ specifies an additional control variable the choice of which is in fact the distribution of the next state and control, whereas the original randomized control $\varphi$ affects only the payoff. This is a team problem as opposed to a control problem because while both controls have the same objective, viz., to maximize a common reward, they are implemented in a non-cooperative manner. This is reminiscent of, e.g., \cite{Hern}, which considers  the cost minimization formulation in which a similar procedure leads  to a zero sum ergodic game. There does not, however, appear to be any corresponding development earlier for the reward maximization problem with a positive reward. While this is completely analogous to the game situation,
we have obtained it without an explicit minorization condition as in \cite{DiMasi}, or the `condition B' of \cite{Anna1}. We have instead conditions \textbf{(A0)} and \textbf{(A1)} which are relatively mild, and compactness of state space, which is not. We are working towards relaxing the latter.\\

An important point to note here is that we have an equivalent problem of maximizing a concave upper semi-continuous function over the convex set $\G_1$. This is in contrast with
 the ergodic \textit{team} problem of maximizing the same function over the \textit{nonconvex} set
\begin{eqnarray*}
\G_2 &:=& \{\eta(dx)\varphi(du|x)\mu(dy|x) : \eta \ \mbox{is invariant under the transition}  \\
&&  \mbox{kernel} \ x \mapsto \mu'(dy|x)\},
\end{eqnarray*}
i.e., where the controls $\varphi, \mu'$ are chosen by the two team members noncooperatively. The latter is what one obtains from the team formulation via log transformation.

\item It  is also worth noting that the entropic penalty implicit in our variational formula also arises in different contexts \cite{Bierkens}, \cite{Guan}, \cite{Todorov}.\\
\end{enumerate}

\section{Examples}			\label{s:examples}

\subsection{Path counting on graphs}

Let $G$ be
a directed graph on a finite vertex set $\SA$ of size $d$,
with edge set denoted by $\mathcal{E}_G$.
Let $M_G$ denote the incidence
matrix of the graph, namely the $d \times d$ nonnegative matrix
$M_G = \left[ m(x,y) \right]$,
with $m(x,y) = 1$ if $(x,y) \in \mathcal{E}_G$, and
$m(x,y) = 0$ otherwise.
Assume that each vertex has at least one out-going edge.
For $n \ge 1$
and $x \in \SA$,
let $N_n(x)$ denote the number of directed paths of length $n$ starting at $x$.
Then the growth rate of the number of directed paths in the graph, namely
\[
\max_{x \in \SA} \lim_{n \to \infty} \frac{1}{n} \log N_n(x)
\]
exists and equals $\log \rho(M_G)$, where
$\rho(M_G)$ is the Perron-Frobenius eigenvalue of $M_G$.

It is also known that this common limit can be written as
\begin{equation}		\label{entropyrate}
\sup_{G-\mbox{\tiny{compatible }} (\Pi, \pi)} - \sum_{x,y} \pi(x) \pi(y|x) \log \pi(y|x)~.
\end{equation}
Here $\Pi$ ranges over $d \times d$ transition probability matrices
that are $G$-compatible for the directed graph $G$, i.e. such that $\pi(y|x) > 0$
implies that $(x,y) \in \mathcal{E}_G$, and $\pi$ ranges over invariant probability distributions
for $\Pi$. Note that this is the largest entropy rate among all stationary Markov
chains whose transition probability matrix is compatible with the graph.

This characterization
of the growth rate of the number of paths in an irreducible graph
is a consequence of the Donsker-Varadhan formula for the Perron-Frobenius eigenvalue of a nonnegative matrix. Let us verify this as a
corollary of Theorem \ref{thm:mainstrong} in the case without controls.
We take the state space in Theorem \ref{thm:mainstrong} to be $\SA$, i.e.
the vertex set of the graph. The control space $U$ is a set consisting
of a single point, which we write as $U = \{u\}$.
Let $p(y|x,u) := \frac{1}{d(x)}$ for $d(x) :=$ the out-degree of $x$
and $(x,y) \in \mathcal{E}_G$, and let

\begin{equation}		\label{rwocontrol}
r(x,u,y) := \begin{cases} \log d(x) &\mbox{ if $(x,y) \in \mathcal{E}_G$}\\
- \infty &\mbox{ otherwise.}
\end{cases}
\end{equation}

Substituting these into the right hand side of (\ref{varstrong}) gives
the expression in (\ref{entropyrate}).

We now bring risk-sensitive
control into this mix of ideas.
Let $U$ be a finite set and suppose now that for each $u \in U$ we are given
a directed graph $G_u$ with vertex set $\SA$.
Assume that each vertex has at least one out-going edge in each $G_u$.
We pose the problem of maximizing
\[
\max_{x \in \SA} \liminf_{n \to \infty} \frac{1}{n} \log \hat{N}_n(x)~,
\]
where now $\hat{N}_n(x)$ is the largest number of directed paths of length $n$ one can
create when starting at $x$ and at each time choosing one of the graphs along
which to move (i.e. one of the control actions) depending on the history
of the states visited so far. More generally, we might allow for a randomized choice
of the graph to be used at each time,
based on the history of the states and the realizations of the control so far, and
ask for the maximum growth
rate of the expectation of the number of directed paths of each length that we
can create in this way.

This problem can be posed in a framework that is amenable to an application of
Theorem \ref{thm:mainstrong}. As in the case without controls, we set
$p(y|x,u) := \frac{1}{d_u(x)}$ for all $(x,y) \in \mathcal{E}_{G_u}$,
where $d_u(x)$ denotes the out-degree of vertex $x$ in $G_u$, and we now set

\begin{equation}		\label{rwithcontrol}
r(x,u,y) := \begin{cases} \log d_u(x) &\mbox{ if $(x,y) \in \mathcal{E}_{G_u}$}\\
- \infty &\mbox{ otherwise.}
\end{cases}
\end{equation}

According to Theorem \ref{thm:mainstrong} this maximum growth rate is given
by
\[
\max_\eta - \sum_{x,u} \tilde{\eta}(x,u)
\sum_{y ~:~ (x,y) \in \mathcal{E}_{G_u}} \eta_2(y|x,u) \log \eta_2(y|x,u))~,
\]
where the maximum is over all $\eta(x,u,y) = \tilde{\eta}(x,u)\eta_2(y|x,u)$
with $\eta_2(y|x,u) > 0$ implying that $(x,y) \in \mathcal{E}_{G_u}$, and
such that
\[
\sum_{(x,u)} \tilde{\eta}(x,u) \eta_2(y|x.u) = \eta_0(x)~,
\]
where, as usual, $\eta_0(x) := \sum_u \tilde{\eta}(x,u)$.
Note that this has following interpretation: among all stationary Markov
chains $((X_n,Z_n), n \ge 0)$ with state space $\SA \times U$ that are compatible
with the
family of graphs in the sense that if a transition from $(x,u)$ to $(y,u^\prime)$
has positive probability then $(x,y) \in \mathcal{E}_{G_u}$, maximize the
conditional entropy of the next state given the current state-entropy pair, i.e.
maximize $H(X_1|X_0,U_0)$.

The interpretation of the growth rate of the number of directed paths of a
given length in a directed graph as an entropy rate has considerable practical
importance in coding theory. Each directed path of length $n$ can be viewed as
an {\em allowed sequence} of length $n$, with coordinates from the state space
$\SA$, and the set of such directed paths is then viewed as a set of
{\em constrained sequences}
\cite[Problem 4.16]{CT}, \cite{LindMar}. The problem of {\em constrained coding} has
been extensively studied. In one version of this problem,
the goal is to come up with algorithms that can
take an infinitely long sequence of symbols from a finite set of size $m$
and produce $\SA$-valued sequences as output
in a one-to-one fashion, and such that the output sequences meet the constraints
defined by the graph, see \cite[Sec. 5.2]{LindMar} for more details.
Naturally, it is not possible to do this if $\log m$ exceeds the growth rate given by (\ref{entropyrate}); finding efficient
algorithms to do this whenever $\log m$ is less than
the growth rate given in (\ref{entropyrate})
was a key early success in this area
\cite{ACH}, \cite{LindMar}.
Investigating the question of constrained
coding up to the maximum possible conditional entropy rate given by the application
of Theorem \ref{thm:mainstrong} to the controlled graph formulation above
would be an interesting challenge.
\color{black}

\subsection{Portfolio optimization}

As another example, we consider the portfolio optimization problem from \cite{Bielecki}, except that we consider the reward maximization framework instead of cost minimization as in the classic work of Cover \cite{Cover}. The model is as follows. The underlying `factor process' $\{X_n\}$ is a discrete time  Markov chain on a finite state space $\Q := \{1, \cdots, m\}$ (say) with an irreducible transition matrix $Q = [[q(j|i)]]$. The control space will be the
simplex
$A := \{a = a_1, \cdots, a_m] \in \R^m : a_i \geq 0 \ \forall i, \ \sum_ia_i \leq 1\}$, with $a_i$ denoting the proportion of wealth invested in the $i$th risky asset. In particular, $1 - \sum_ia_i$ is then the proportion invested in the risk-less bank account. We denote by $\{\pi_n\}$ the $A$-valued control sequence, representing the trading strategy, i.e., $\pi_{n,i}$ will be the proportion of wealth invested in the $i$th risky asset at time $n$.  $\{W_n\}$ is the process of $m$-dimensional vectors of price relatives such that  $W_{n+1}$ is conditionally independent of $X_i, i < n; W_i, \pi_i, i \leq n$,   given $(X_n, X_{n+1})$ and its conditional law given the latter is specified by a kernel $\nu(x, y, dw) : \Q\times\Q \mapsto \R^m$ with support in the interior of the positive cone in $\R^m$. Let $e^r, r > 0,$ denote the per period multiplier of wealth invested in the bank account (thus $e^r - 1$ is the interest rate). Let $\1$ denote the constant vector of all $1$'s. The evolution of the wealth process $\{V_t\}$ is given by
\begin{displaymath}
V_{n+1} = V_n\left[e^r + \langle\pi_n, W_{n+1} - e^r\1\rangle\right]~,
\end{displaymath}
where $V_0 := 1$.
The objective is to maximize the
risk-adjusted growth rate of wealth

\begin{equation}
\liminf_{n\uparrow \infty}\frac{1}{n}\log E\left[e^{- \frac{\theta}{2}\log V_n}\right]. \label{finrew}
\end{equation}

Here
$\theta$
is the risk sensitivity parameter.
The control sequence $\{\pi_n\}$ is assumed to be adapted to the
factor process $\{X_n\}$ and the
controls, i.e. the distribution of
$\pi_n$ is chosen as a function of
$(X_0, \ldots, X_n, \pi_0, \ldots, \pi_{n-1})$.\\

It is useful to constrast the objecive we consider with that
considered in \cite{Bielecki} of
maximizing, for $\theta > 0$, the quantity
\begin{equation}
\liminf_{n\uparrow \infty}-\frac{2}{\theta}\frac{1}{n}\log E\left[e^{-\frac{\theta}{2}\log V_n}\right]~. \label{finrewcost}
\end{equation}
In \cite{Bielecki}, this problem is considered by writing the objective in (\ref{finrewcost})
as
\[
- \limsup_{n\uparrow \infty}
\frac{2}{\theta}\frac{1}{n}\log E\left[e^{-\frac{\theta}{2}\log V_n}\right]~,
\]
and then studying the risk-sensitive cost minimization problem corresponding to the objective
\[
\limsup_{n\uparrow \infty}
\frac{2}{\theta}\frac{1}{n}\log E\left[e^{-\frac{\theta}{2}\log V_n}\right]~.
\]
That positive $\theta$ indicates risk aversion in (\ref{finrewcost})
is argued, see \cite[Eqn. (2.1)]{BPS},
by writing the Taylor's series expansion, for small $\theta$,
\[
-\frac{2}{\theta} \log E\left[e^{-\frac{\theta}{2}\log V_n}\right]
= E[ \log V_n] - \frac{\theta}{4} \var(\log V_n) + O(\theta^2)~.
\]

By constrast, our formulation is able to handle both the case of risk-aversion and risk-seeking.
The Taylor's series expansion
\[
\log E\left[e^{-\frac{\theta}{2}\log V_n}\right] = - \frac{\theta}{2} E[ \log V_n]
+ \frac{\theta^2}{8} \var(\log V_n) + o(\theta^2)
\]
indicates that if the objective in (\ref{finrew}) is multiplied by $- \frac{2}{\theta}$, then
it corresponds to risk-aversion for positive $\theta$ and to risk-seeking for
negative $\theta$.


Keeping in mind that $e^r + \langle a, W - e^r\textbf{1} \rangle > 0$  under our assumption on the support of $\nu(x,y,dz)$, define

\begin{eqnarray*}
\mu(x,a,y) &:=& \int e^{- \frac{\theta}{2}\log\left[e^r + \langle a, w - e^r\rangle\right]}\nu(x, y, dw), \\
&& (\mbox{assumed to be} \ < \infty) \\
r(x, a) &:=& - \frac{2}{\theta}\log\left(\sum_yq(y|x)\mu(x,a,y)\right), \\
p(y|x, a) &:=& \frac{q(y|x)\mu(x,a,y)}{\sum_{y'}q(y'|x)\mu(x,a,y')}.
\end{eqnarray*}

One can show that for all $n \ge 1$ and all admissible controls, we have
\begin{displaymath}
\frac{1}{n}\log E\left[e^{- \frac{\theta}{2}\log V_n}\right] = \frac{1}{n}\log \tE\left[e^{- \frac{\theta}{2}\sum_{m=0}^{n-1}r(X_m, \pi_m)}\right]~,
\end{displaymath}
where $\tE$ is the expectation with respect to the law
\begin{eqnarray*}
\lefteqn{p(x_0) \phi_0(da_0|x_0) p(x_1|x_0,a_0) \phi_1(da_1|x_0, a_0, x_1) \ldots} \\
&\times&\phi_{n-1}(da_{n-1}|(x_i, a_i, 0 \le i \le n-2),x_{n-1})~,
\end{eqnarray*}
where $p(x_0)$ is the initial distribution of $X_0$, the admissible controls are determined
by the kernels $\phi_0( \cdot|\cdot), \ldots, \phi_{n-1}(\cdot|\cdot)$, and the salient point
is that the transition kernel for the evolution of the factor process under this change of measure
is given by the kernel $p(\cdot|\cdot, \cdot)$ defined above.
To see this, first observe that $W_1, \ldots, W_n$ are conditionally independent and
identically distributed given $(X_i, \pi_i, 0 \le i \le n)$. Hence
\begin{eqnarray*}
E\left[e^{- \frac{\theta}{2}\log V_n}| X_i, \pi_i, 0 \le i \le n \right]
&=& E\left[\prod_{m=0}^{n-1} e^{- \frac{\theta}{2}\log \frac{V_{m+1}}{V_m}}
| X_i, \pi_i, 0 \le i \le n \right]\\
&=& \prod_{m=0}^{n-1}E\left[e^{- \frac{\theta}{2}\log \frac{V_{m+1}}{V_m}}
| X_i, \pi_i, 0 \le i \le n \right]\\
&=& \prod_{m=0}^{n-1} \mu(X_m, \pi_m, X_{m+1})~,
\end{eqnarray*}
so we have
\[
 E\left[e^{- \frac{\theta}{2}\log V_n}\right]
 = E\left[\prod_{m=0}^{n-1} \mu(X_m, \pi_m, X_{m+1})\right]~.
\]
For an admissible control strategy, we can write this as
\begin{eqnarray*}
&& \sum_{x_0, \ldots, x_n} \int_{a_0} \ldots \int_{a_{n-1}} p(x_0) \prod_{m=0}^{n-1}
\mu(x_m,a_m,x_{m+1}) q(x_{m+1}|x_m)\\
&&~~~~~~~~~~~~~~~~~~~~~~~~~~~~~~~~~~~~~~~~~~ \phi_m(d a_m|(x_i, a_i, 0 \le i \le m-1),x_m)~,
\end{eqnarray*}
which is the same as
\begin{eqnarray*}
&& \sum_{x_0, \ldots, x_n} \int_{a_0} \ldots \int_{a_{n-1}} p(x_0) \prod_{m=0}^{n-1}
e^{- \frac{\theta}{2} r(x_m, a_m)} p(x_{m+1}|x_m, a_m)\\
&&~~~~~~~~~~~~~~~~~~~~~~~~~~~~~~~~~~~~~~~~~~ \phi_m(d a_m|(x_i, a_i, 0 \le i \le m-1),x_m)~,
\end{eqnarray*}
which equals $\tE\left[e^{- \frac{\theta}{2}\sum_{m=0}^{n-1}r(X_m, \pi_m)}\right]$.



Hence the problem of maximizing (\ref{finrew}) is equivalent to the risk-sensitive control problem for a controlled Markov chain on $\Q$ with action space $A$ and controlled transition probabilities $p(y|x, a), x,y \in \Q, a \in A$, the objective being to maximize the reward

\begin{displaymath}
\lambda := \sup_{x_0} \sup \liminf_{n\uparrow\infty} \frac{1}{n}\log E\left[e^{- \frac{\theta}{2}\sum_{m=0}^{n-1}r(X_m, \pi_m)} | X_0 = x_0 \right]~.
\end{displaymath}
where the second supremum is over admissible controls.

The optimal growth rate for the wealth is then given by

\begin{displaymath}
\lambda = \max_{\eta \in \mathcal{G}}\Big(\sum_x\int_A \tilde{\eta}(x,da)(- \frac{\theta}{2}r(x,a) - \sum_y\int_A \eta_2(y|x,a)\log\left(\frac{\eta_2(y|x,a)}{p(y|x,a)}\right))\Big)
\end{displaymath}

where

\begin{eqnarray*}
\mathcal{G} &:=& \Big\{\eta(x,da,y) \in \PA(\Q\times A \times \Q) : \eta(x,da,y)
= \tilde{\eta}(x,da) \eta_2(y|x, a) \\
 && = \eta_0(x)\eta_1(da|x)\eta_2(y|x, a) \ \mbox{such that} \ \eta_0 \ \mbox{is stationary under} \\
 && \mbox{the transition matrix} \ \left[\left[\int\eta_1(da|x)\eta_2(y|x,a)\right]\right]\Big\}.
 \end{eqnarray*}
 

 In order to justify this, we need to verify that the conditions
 \textbf{(A0)} and \textbf{(A1)} are satisfied. Here $\Q$ plays the role of
 $\SA$, $A$ plays the role of $U$, and $-\frac{\theta}{2}r(x,a)$ plays the
 role of $r(x,u,y)$ in the general theory. The validity of \textbf{(A0)} follows
 from the continuity of the logarithm function. The validity of \textbf{(A1)}
 follows from the continuity of the logarithm function, the fact that $\Q$ is finite,
 and because $\sum_{y'} q(y'|x) \mu(x,a,y)$ is strictly positive for all $(x,a)$.

If we discretize $A$, this is a finite dimensional concave maximization problem eminently amenable to standard nonlinear programming tools.

\subsection{Minimizing exit rate from a domain}

 Consider a  set of controlled stochastic matrices on a finite state space $S = \{1, \cdots, s\}$
denoted by $P_u = [[p(j|i, u)]]_{i,j \in S}$.
Here $u$ is the control parameter taking values in $A$,
 where $A$ is a compact metric action space. We assume that $u \mapsto P_u$ is continuous and $P_u$ is irreducible for all $u$.  Let $S_0 \subset S$ be
a nonempty proper subset of $S$
and let $S_1 := S_0^c$ denote its complement.
 Let $\check{P}_u$ denote the restriction of $P_u$ to $S_1$ and
for a sequence of random variables $\{X_n\}$ with values in $S$,
define $\tau := \inf\{n \geq 0 : X_n \in S_0\}$.

We are interested in determining
\[
\lambda := \sup_{i \in S_1} \sup \liminf_{n\uparrow\infty}\frac{1}{n}\log P(\tau > n)~,
\]
where the second supremum is over all admissible controls, and the law of $\tau$ is
determined by the control strategy.
Namely, we are interested in the problem of finding the slowest exit rate from $S_1$ over admissible
control strategies.

Write $\check{P}_u = D_uQ_u$ where $D_u$ is a diagonal matrix with its $i$th diagonal entry $d(i,u) := \sum_{j \in S_1}p(j|i,u)$ and $Q_u := [[q(j|i,u)]]$ is
a stochastic matrix on $S_1$
given by $q(j|i,u) := d(i,u)^{-1}p(j|i,u)$, where we will also assume that
$d(i,u) > 0$ for all $i \in S_1$ and $u \in A$.
It can be checked that for any
admissible control strategy and $i \in S_1$, we have
\[
P(\tau > n) = E\left[e^{\sum_{m=0}^{n-1}\log(d(X_m, U_m))}\right]~,
\]
where $U_m$ denotes the choice of control at time $m$, and
$\{X_n\}$ is the $S_1$-valued Markov chain, having the transition probability matrix
$Q_{U_m}$ at time $m$. Therefore, with the choices $\SA := S_1$, $U := A$,
and $r(i,u,j) := \log d(i,u)$, the problem is amenable to our general theory.


 Disintegrate a typical element
$\eta \in \PA(S_1 \times A\times S_1)$ as
$\eta_0(i) \eta_1(du|i)\eta_2(j | i, u)$, and write
$\tilde{\eta}(i, du)$ for $\eta_0(i) \eta_1(du|i)$.

 Then our results show that

 \begin{eqnarray*}
&&\lambda = \max_{\eta \in \mathcal{G}}\Big(\sum_{i,j \in S_1}\int_{A}\eta(i, du, j)\log(d(i,u)) - \\
&&~~~~~~~~~~~~\sum_{i \in S_1}\int_{A}\tilde{\eta}(i, du)D(\eta_2(j| i, u)||q(j|i,u))\Big)~,
\end{eqnarray*}
where $\mathcal{G}$ denotes the set of $\eta \in \PA(S_1 \times A\times S_1)$ for which
$\eta_0$ is invariant under the transition kernel $\int_A \eta_1(du|i)\eta_2(j | i, u)$.
To verify this, we need to check the validity of the conditions
\textbf{(A0)} and \textbf{(A1)}. The former is a consequence of the assumed continuity of
$u \mapsto P_u$. The latter is a consequence of the fact that $S_1$ is finite and that
$u \mapsto Q_u$ is continuous, which in turn follows from the assumed continuity of
$u \mapsto P_u$ and the assumption that $d(i,u) > 0$ for all $i \in S_1$ and
$u \in A$.

\section{Concluding remarks}		\label{s:concremarks}

We considered the problem of maximizing the growth rate of reward in the standard risk-sensitive formulation for a controlled Markov chain on a compact metric state space, with a compact metric action space. We took a non-standard approach to this problem via  a nonlinear version of the Krein-Rutman theorem
to obtain a variational formulation for the optimal reward. This leads to an occupation measure based concave maximization formulation of the control problem.\\

 The approach holds promise for possible use of convex optimization techniques for approximate solution of the risk-sensitive reward maximization problem, in a manner analogous to what abstract linear programming does for the classical additive reward problems (such as discounted or ergodic rewards, see, e.g., \cite{HeLa}). We achieved this with rather few technical conditions except for the compactness of the state and action spaces. It remains a major challenge to extend this approach to noncompact state and action spaces.

\end{document}